\documentclass[12pt,reqno]{amsart}

\usepackage{amsfonts,amsmath,amsthm}
\usepackage{amssymb,epsfig}
\usepackage{enumerate} 

\usepackage{color} %color
\definecolor{vert}{rgb}{0,0.6,0}

\usepackage[text={425pt,650pt},centering]{geometry}
\usepackage{graphicx}
\usepackage{epsfig}
\usepackage{tikz,esvect}
\usetikzlibrary{3d,calc,shapes}

\usepackage{caption}
\usepackage{color} %color
\definecolor{vert}{rgb}{0,0.6,0}

\numberwithin{figure}{section}

\theoremstyle{plain}
\newtheorem{thm}{Theorem}[section]

\newtheorem{lem}[thm]{Lemma}

\theoremstyle{remark}
\newtheorem{rem}{\bf{Remark}}
\numberwithin{equation}{section}

\newcommand{\N}{\mathbb{N}}

\newcommand{\R}{\mathbb{R}}

\newcommand{\T}{\mathbb{T}}
\newcommand{\Z}{\mathbb{Z}}

\newcommand{\cA}{\mathcal{A}}
\newcommand{\cB}{\mathcal{B}}

\newcommand{\cS}{\mathcal{S}}

\newcommand{\AC}{{\rm AC\,}}

\newcommand{\BUC}{{\rm BUC\,}}

\newcommand{\Lip}{{\rm Lip\,}}

\newcommand{\al}{\alpha}
\newcommand{\gam}{\gamma}

\newcommand{\ep}{\varepsilon}

\newcommand{\ol}{\overline}

\begin{document}

\title[Optimal rate of convergence]
{Optimal convergence rate for periodic homogenization of convex Hamilton-Jacobi equations}

\author[H. V. TRAN, Y. YU]
{Hung V. Tran, Yifeng Yu}

\thanks{
The work of HT is partially supported by NSF CAREER grant DMS-1843320, a Simons Fellowship, and a Vilas Faculty Early-Career Investigator Award.
The work of YY is partially supported by NSF grant 2000191.
}

\address[H. V. Tran]
{
Department of Mathematics, 
University of Wisconsin Madison, Van Vleck Hall, 480 Lincoln Drive, Madison, Wisconsin 53706, USA}
\email{hung@math.wisc.edu}

\address[Y. Yu]
{
Department of Mathematics, 
University of California at Irvine, 
California 92697, USA}
\email{yyu1@math.uci.edu}

%\date{\today}

\maketitle

\begin{abstract}
In this paper,  we show that the  rate of convergence  in periodic homogenization of convex Hamilton-Jacobi equations is  always $O(\ep)$, which is optimal.   
This is a natural extension of a result concerning stable norms  in metric geometry  \cite{Burago}  that is essentially equivalent to the homogenization of  convex static Hamilton-Jacobi equations.  
Another extremely interesting question in this direction is whether the $O(\ep)$ rate holds in the nonconvex setting.   
We present  a special nonconvex example with  $O(\ep)$ convergence rate, which relies on identifying the shape of the effective Hamiltonian and game theory interpretation formulas.

\end{abstract}

%%%%%%%%%%%%%%%%%%%%%%%%%%%%%%%%%%%%%%%%%%%%%%%%%%%%%%%%%%%%%%%%%%%%%%%%%%%%%%%%%%%%%%%%

\section{Introduction}

We give a brief description of the periodic homogenization of Hamilton-Jacobi equations.
For each $\ep>0$, let $u^\ep \in C(\R^n \times [0,\infty))$ be the viscosity solution to 
\begin{equation}\label{C-ep}
\begin{cases}
 u_t^\ep+H\left(\frac{x}{\ep},Du^\ep\right)=0 \qquad &\text{in} \ \R^n \times (0,\infty),\\
u^\ep(x,0)=g(x) \qquad &\text{on} \ \R^n.
\end{cases} 
\end{equation}
Here, the Hamiltonian $H=H(y,p):\R^n \times \R^n \to \R$ is a  given continuous function  satisfying
\begin{itemize}
\item[(H1)] for $p\in \R^n$, $y \mapsto H(y,p)$ is $\Z^n$-periodic;
\item[(H2)] $H$ is coercive in $p$, that is, uniformly for $y\in \T^n:=\R^n/\Z^n$,
\[
\lim_{|p| \to \infty} H(y,p) = +\infty;
\]
\item[(H3)] for $y\in \T^n$, $p \mapsto H(y,p)$ is convex.

\end{itemize}
The initial data $g\in \BUC(\R^n) \cap \Lip(\R^n)$, where  $\BUC(\R^n)$ is the set of bounded, uniformly continuous functions on $\R^n$.

Under assumptions (H1)--(H2), $u^\ep$ converges to $u$ locally uniformly on $\R^n \times [0,\infty)$ as $\ep \to 0$, and $u$ solves the effective equation
\begin{equation}\label{C}
\begin{cases}
 u_t+\ol{H}\left(Du\right)=0 \qquad &\text{in} \ \R^n \times (0,\infty),\\
u(x,0)=g(x) \qquad &\text{on} \ \R^n.
\end{cases} 
\end{equation}
See \cite{LPV,Ev1,Tran}.
The effective Hamiltonian $\ol{H} \in C(\R^n)$ depends nonlinearly on $H$, and is determined by the cell (ergodic) problems.  
For each $p \in \R^n$, there exists a unique constant $\ol H(p)\in \R$ such that the following cell problem has a continuous viscosity solution
\begin{equation}\label{E-p}
H(y,p+Dv)=\ol H(p)  \qquad \text{in} \ \T^n.
\end{equation}
Note that $v=v(y,p)$ is not unique even up to additive constants in general.

Our goal in this paper is to obtain rate of convergence of $u^\ep$ to $u$ in $L^\infty$, that is, an optimal bound for $\|u^\ep-u\|_{L^\infty(\R^n \times [0,\infty))}$  as $\ep \to 0+$.  
Here is our first main result.

\begin{thm} \label{thm:main}
Assume {\rm (H1)--(H3)} and $g \in \BUC(\R^n) \cap \Lip(\R^n)$.
For $\ep>0$, let $u^\ep$ be the viscosity solution to \eqref{C-ep}.
Let $u$ be the viscosity solution to \eqref{C}.
Then, there exists $C>0$ depending only on $H$, $\|Dg\|_{L^\infty(\R^n)}$, and $n$ such that
\[
\|u^\ep-u\|_{L^\infty(\R^n \times [0,\infty))} \leq C \ep.
\]
\end{thm}

\subsection{Brief history} 
We now give a minimalistic review of the PDE literature  playing major roles in finding the convergence rate.

\begin{enumerate}

\item  For the general nonconvex setting, the best known convergence rate is $O(\ep^{1/3})$ obtained by Capuzzo-Dolcetta and Ishii \cite{CDI}.

\item For convex Hamilton-Jacobi equation, by using weak KAM methods,  the lower bound $u^\ep - u \geq -C\ep$  was proved by Mitake, Tran, and Yu \cite{MTY} for all dimensions. 
When $n=2$,  for homogeneous Hamiltonian (e.g., those  physically relevant ones from front propagations), the upper bound  $u^\ep - u \leq C\ep$ was also derived via the classical Aubry-Mather theory, which  however heavily relies on the two dimensional topology.     After \cite{MTY},  the major open problem was whether  the upper bound $u^\ep - u \leq C\ep$ always holds in higher dimensions ($n\geq 3$), which was completely unclear although some conditional higher dimensional results in \cite{MTY} sort of  imply that the upper bound should hold in ``generic" situations.

\item  Recently,  Cooperman  has discovered that closely related convergence rate results have been established in the context of first passage percolation (FPP) decades ago (\cite{A1,A2}).  
By adjusting the methods there, he was able to obtain in \cite{Cooperman}  a near optimal convergence rate  $|u^\ep(x,t) - u(x,t)| \leq C\ep \log (C+\ep^{-1}t)$ when $n\geq 3$, which was a quite surprising result for  people in the PDE community. 

\item  Later on,  while studying finer properties of  effective fronts,  the authors discovered an old result concerning the convergence rate of stable norms in metric geometry by Burago \cite{Burago}, that is basically equivalent to the optimal convergence rate $O(\ep)$ in the homogenization of static Hamilton-Jacobi equations.  
This can be easily extended to the Cauchy problem via standard adjustments, which concludes the study of this whole program in the convex setting.

\end{enumerate}

Interestingly,   the word ``homogenization" rarely appeared in the stable norm literature  in spite of  the essential equivalence. 
Similarly, related important results from stable norm literature are hidden from people in the PDE community where the focus is more on the wellposedness theory.    
The research focus in the relevant geometry and dynamical system communities is mainly on characterizing finer properties of stable norms (see \cite{Burago, TY2022} and reference therein), which might be part of the reasons that quantitative results like the convergence rate went unnoticed by PDE community.  
From a technical point of view,   it is probably fair to say that the optimal convergence problem in the convex setting has already been solved in \cite{Burago} in an equivalent form.   
For the reader's convenience, we give a brief review of stable norms and corresponding results and proofs in \cite{Burago} where connections are  more obvious than those in FPP literature  (\cite{A1,A2}). 

\subsection{Connection with stable norms} 
For the clarity of the presentation, let us focus on a relatively  simple form of stable norms.  
For $a\in C(\Bbb T^n, (0,\infty))$,  we define a corresponding periodic Riemannian metric on $\R^n$ as
\[
g =\frac{1}{a(x)}\sum_{i=1}^{n}dx_{i}^2.
\]  
Let $d_a(\cdot,\cdot)$ denote the distance function induced by this metric.
The \emph{stable norm} associated with $g$ (or $a$) is defined as
\begin{equation}\label{x-a}
\|x\|_{a}=\lim_{\lambda\to \infty}{d_a(0,\lambda x)\over \lambda} \qquad \text{ for } x\in \R^n.
\end{equation}
Stable norms and their properties in various settings have been extensively studied in the communities of dynamical systems and geometry.  
See \cite{Burago} for more background. 
The limit \eqref{x-a}  can also be viewed as a homogenization problem of the following static Hamilton-Jacobi equation
\[
\begin{cases}
a\left({x\over \ep}\right)|Dw^{\ep}|=1 \qquad \text{ for $x\in \R^n\backslash\{0\}$,}\\
w^{\ep}(0)=0.
\end{cases}
\]
Here, $w^\ep$ is the maximal viscosity solution to the above.
By the optimal control formula, for $x\in \R^n$,
\[
w^\ep(x)=\ep d_a\left(0,\frac{x}{\ep}\right).
\]
As $\ep \to 0$, $w^\ep \to w =\|\cdot\|_a$ locally uniformly on $\R^n$, and $w=\|\cdot\|_a$ is the maximal viscosity solution of 
\[
\begin{cases}
\ol H_a(Dw)=1 \qquad \text{ for $x\in \R^n\backslash\{0\}$,}\\
w(0)=0.
\end{cases}
\]
Here, $\ol H_a$ is the effective Hamiltonian of $a(y)|p|$.
 It was proved in  \cite{Burago} that a periodic Riemannian distance is not much different from the Euclidean distance, i.e., 
\begin{equation}\label{optimal-rate}
\left|\lambda \|x\|_a-d_a(0,\lambda x)\right|\leq C \qquad \text{ for all } x\in \R^n, \lambda>0,
\end{equation}
where $C>0$ is a universal  constant.  
Note that in PDE language  this  is actually equivalent to  the optimal convergence rate
\[
\|w^{\ep}-w\|_{L^\infty(\R^n)}\leq C\ep.
\]
As an immediate corollary,  for homogenization of the static problem in a bounded open set $U \subset \R^n$, 
\[
\begin{cases}
a\left({x\over \ep}\right)|Du^{\ep}|=1 \qquad &\text{ in } U,\\
u^{\ep}(x)=g(x)  \qquad &\text{ on }  \partial U,
\end{cases}
\]
and the effective equation
\[
\begin{cases}
\ol H_a(Du)=1 \qquad &\text{ in } U,\\
u(x)=g(x)   \qquad &\text{ on }  \partial U,
\end{cases}
\]
we have the convergence rate 
\[
|u^\ep(x)- u(x)|\leq C\ep.
\]
In fact, without loss of generality, it suffices to assume that $0\in U$ and establish the above for $x=0\in U$. 
The convergence rate follows directly from  the optimal control formulation for static equations (\cite{L})
\[
u^{\ep}(0)=\min_{x\in \partial U}\left\{g(y)+w^{\epsilon}(y)\right\}.
\]
The connection between  Hamilton-Jacobi equations and geometry is well-known to experts. The general static Hamilton-Jacobi equation is corresponding to Finsler metric. See \cite{LN} for instance. 

The proof of \eqref{optimal-rate} is divided into two steps in  \cite{Burago}. 

\medskip

\noindent {\bf Step 1.}  It is quite standard  to have the subadditive property  (\cite[Lemma 3]{Burago})
\begin{equation}\label{lowerdouble}
d_a(0, 2x)\leq 2d_a(0,x)+C \qquad \text{ for } x\in \R^n,
\end{equation}
which leads to the lower  bound $d_a(0,\lambda x)\geq \lambda \|x\|_a-C$ by a standard   iteration (\cite[Lemma 1]{Burago}). 

\medskip
 
\noindent {\bf Step 2.}  By using a beautiful curve cutting lemma (\cite[Lemma 4]{Burago}), one obtains the superadditive property
\begin{equation}\label{upperdouble}
d_a(0, 2x)\geq 2d_a(0,x)-C. 
\end{equation}
Then, a standard   iteration gives the lower bound $d_a(0,\lambda x)\leq \lambda \|x\|_a+C$. 

\medskip

The subadditive inequality \eqref{lowerdouble} is rather natural.  
Note that the corresponding lower bound can be obtained by other methods (see e.g., \cite{MTY}).   
The main difficulty is to obtain the superadditive property \eqref{upperdouble}, which is very surprising. 
In fact, a minimizing geodesic (or an action minimizing curve) connecting $0$ and $2x$ does not need to pass through $x$ or within a uniformly bounded distance from $x$ in general.  
When $n=2$,   orbits in the Mather set are known to be within finite distance of a straight line \cite{B} (bounded deviation), which was employed in \cite{MTY} to obtain $O(\ep)$ (essentially \eqref{lowerdouble}) in two dimensions. 
However, the property of bounded deviation  is  no longer true when $n\geq 3$. 

    In \cite{Burago},  \eqref{upperdouble} was proved by a beautiful  topological lemma on decompositions of continuous curves (Lemma \ref{lem:top}) that has nothing to do with dynamics.  As an immediate corollary,  every  minimizing geodesic connecting $0$ and $2x$ can be cut into different pieces and then reassembled into a new near-minimizing geodesic that passes through $x$.  For the reader's convenience,  we present a proof of Lemma \ref{lem:top} in Section \ref{sec:prelim}.  The method in \cite{Burago} is quite robust.  Due to the similar optimal control formulation and the well-known metric form of the Lagrangian,  the extension  to the Cauchy problem is rather straightforward for experts. 

As it was aforementioned,   similar results (with an extra $\log$ correction term) and proofs have also appeared in \cite{A1,A2} in the study of the limiting shape arising from FPP.  
In fact, the setting of  FPP corresponds to the case where $a(x)$ is an i.i.d  distribution on the integer lattice $\Z^n$. 
The limit shape in FPP corresponds to the unit ball of the stable norm (equivalently,  the effective front in the context of homogenization).  However, unlike  \cite{Burago},   ideas and methods in \cite{A1,A2} were not evidently close to our context. 
In this aspect,   \cite{Cooperman} plays a significant role in finding and organizing  them   into  forms applicable to Hamilton-Jacobi equations in the periodic environment. 

\subsection{A nonconvex example}  
The main remaining open problem in the quantitative theory is to find optimal convergence rates for periodic homogenization of nonconvex Hamilton-Jacobi equations, where all methods and ideas from stable norms and FPP cease to work.  
Unlike the case of convex Hamiltonians,  there are different types of nonconvexity. 
The initial data might also play some roles since for linear initial data $g(x)=p\cdot x$ for a given $p\in \R^n$, the convergence rate is always $O(\ep)$ for all $H$.   
The first step towards resolving this open problem is to find nontrivial interesting examples where $O(\ep)$ rate holds or fails.  
As a starting point, we present in this paper a nonconvex example where $O(\ep)$ holds.

\begin{thm} \label{thm:main2}
Assume
\[
H(y,p)=\max\{|p|-1, 1-|p|\}+V(y),
\]
where $V\in C(\T^n)$ with $\max_{\T^n}V-\min_{\T^n}V\geq 1$.
Assume $g\in \Lip(\R^n)$. 
Then, there exists $C>0$ depending only on $\|Dg\|_{L^\infty(\R^n)}$ and $V$ such that
\[
\|u^\ep-u\|_{L^\infty(\R^n \times [0,\infty))} \leq C \ep.
\]

\end{thm}

Theorem \ref{thm:main2} holds  for more general Hamiltonians $H$.   
We choose the specific form $H(y,p)=\max\{|p|-1, 1-|p|\}+V(y)$ to demonstrate the main ideas instead of delving into technicalities.

The proof of Theorem \ref{thm:main2} relies on two main points. 
First,  we use  \cite{ATY1} to identify the structure of $\ol H$. 
In particular, for the large oscillation case ${\rm osc}(V) =\max_{\T^n}V-\min_{\T^n}V\geq 1$, the effective Hamiltonian is convex, which is a special example of the so-called ``quasiconvexification"  phenomenon (see \cite{ATY2}).   
Then, we apply suitable game theory interpretations and comparison principles to link this nonconvex case to a proper convex case.  
It is not completely clear to us whether the result still holds for the small oscillation case ${\rm osc}(V) =\max_{\T^n}V-\min_{\T^n}V<1$. 
See Remark \ref{smallosc} for more discussions.

\subsection*{Organization of the paper}
The paper is organized as follows.
Some preliminary results are given in Section \ref{sec:prelim}.
We then give the proof of Theorem \ref{thm:main} and the proof of Theorem \ref{thm:main2} in Section \ref{sec:proof} and Section \ref{sec:proof2}, respectively.

%%%%%%%%%%%%%%%%%%%%%%%%%%%%%%%%%%%%%%%%%%%

\section{Preliminaries} \label{sec:prelim}
We assume the setting of Theorem \ref{thm:main}.

\subsection{Simplifications} 
We have the following simplifications, which are similar to those in \cite{MTY}.
We use the comparison principle to obtain that
\[
\|Du^{\ep}\|_{L^{\infty}(\R^n\times [0,\infty))}\leq C_0.
\]
Here, $C_0>0$ is a constant depending only on $H$ and $\|Dg\|_{L^\infty(\R^n)}$.  
See \cite{Tran}.
Thus, the values of $H(y,p)$ for $|p|>C_0$ are irrelevant.  
Hence, by modifying $H(y,p)$ for $|p|>2C_0+1$ if needed, we assume further that $H \in C^2(\T^n \times \R^n)$, $H$ grows quadratically in $p$, that is,
\begin{equation}\label{quad-growth}
 {1\over 2}|p|^2-K_0\leq H(y,p)\leq {1\over 2} |p|^2+K_0 \quad \text{ for all } (y,p) \in \T^n \times \R^n,
\end{equation}
for some $K_0 >1$.  
Let $L(y,q)$  be the Lagrangian (Legendre transform) of the Hamiltonians $H(y,p)$. 
It is clear that
\begin{equation}\label{quad-growth-L}
 {1\over 2}|q|^2-K_0\leq L(y,q) \leq {1\over 2} |q|^2+K_0 \quad \text{ for all } (y,q) \in \T^n \times \R^n.
\end{equation}

\smallskip

For $(x,t)\in \R^n \times (0,\infty)$, the optimal control formula for the solution to \eqref{C-ep} gives
\begin{equation}\label{oc}
u^{\ep}(x,t)=\inf_{\substack{ \ep \eta(0)=x \\ \eta \in \AC([-\ep^{-1}t,0])}}\left\{g\left(\ep\eta \left(-\ep^{-1}t\right)\right)+\ep\int_{-\ep^{-1}t}^{0}L(\eta(s),\dot \eta(s))\,ds\right\}. 
\end{equation}
Here, $ \AC([-\ep^{-1}t,0])$ is the space of absolutely continuous curves from $[-\ep^{-1}t,0]$ to $\R^n$.

\subsection{A topological lemma}
The following is a topological lemma from \cite{Burago}. 
The proof below is attributed to both Burago and Perelman. 

\begin{lem}\label{lem:top}
Let $m\in \N$ and $\xi:[0,1] \to \R^m$ be a continuous path.
Then, there is a collection of disjoint intervals $\{[a_i,b_i]\}_{1 \leq i \leq k} \subset [0,1]$ with $k \leq \frac{m+1}{2}$ such that
\[
\sum_{i=1}^k (\xi(b_i)-\xi(a_i)) = \frac{\xi(1)-\xi(0)}{2}.
\]
\end{lem}
This is basically a generalized version of the intermediate value theorem in multi dimensions.
We include the proof here for completeness.

\begin{proof}
Consider the unit sphere $\cS^m$ in $\R^{m+1}$.
For $x=(x_1,\ldots,x_{m+1}) \in \cS^m$, we have $\sum_{i=1}^{m+1} x_i^2=1$.
We now define a map $f:\cS^m \to \R^{m}$ as following.

\smallskip

For each $x\in \cS^m$, take a partition $0=t_0 \leq t_1 \leq \ldots  \leq t_{m+1}=1$ such that
\[
t_i - t_{i-1} = x_i^2 \quad \text{ for } 1\leq i \leq m+1.
\]
Denote by
\[
f(x) = \sum_{i=1}^{m+1} {\rm sign}(x_i) (\xi(t_i) - \xi(t_{i-1})).
\]
Note that if $x_i=0$, then $t_{i-1}=t_i$.
It is clear that $f \in C(\cS^m, \R^{m})$ and $f$ is odd, that is,
\[
f(x)=-f(-x) \quad \text{ for } x\in \cS^m.
\]
Thus, by the Borsuk--Ulam theorem, there exists $x\in \cS^m$ such that 
\[
f(x)=f(-x) \quad \Longrightarrow \quad f(x)=0.
\]
Without loss of generality, we assume that $x$ has at most $\frac{m+1}{2}$ positive coordinates.
Then, the collection of disjoint intervals $[t_{i-1},t_i]$ with $x_i>0$ is exactly what we need.
The proof is complete.
\end{proof}

\section{Proof of Theorem \ref{thm:main}}\label{sec:proof}

We assume the setting of Theorem \ref{thm:main}.
By the simplifications above, we also assume \eqref{quad-growth}--\eqref{quad-growth-L}.
For $x,y\in \R^n$ and $t>0$, denote by
\[
m(t,x,y) = \inf \left\{ \int_0^t L(\eta(s),\dot \eta(s))\,ds\,:\, \eta \in \AC([0,t],\R^n), \eta(0)=x, \eta(t)=y \right\}.
\]
Here, $m(t,x,y)$ is the minimum cost to travel from $x$ to $y$ in a given time $t>0$.
We say that $m(t,x,y)$ is the metric distance from $x$ to $y$ in time $t$.
The homogenized (large time average) metric is 
\[
\ol m(t,x,y) = \lim_{k \to \infty} \frac{1}{k} m(kt,kx,ky).
\]
In fact, 
\[
\ol m(t,x,y) = t\ol L\left( \frac{y-x}{t}\right),
\]
where $\ol L$ is the Lagrangian (Legendre transform) of the effective Hamiltonian $H$.
In particular, for $s>0$,
\[
\ol m(st,sx,sy) = st\ol L\left( \frac{y-x}{t}\right) = s \ol m(t,x,y).
\]

Some basic properties of $m$ are collected in the following lemma.

\begin{lem}\label{lem:basic-m}
Assume \eqref{quad-growth}--\eqref{quad-growth-L}.
The following properties hold, where $C>0$ is a universal constant depending only on $L$ and $n$.
\begin{itemize}
\item[(a)] $m$ is subadditive, that is, for $x,y,z\in \R^n$ and $t,s>0$,
\[
m(t,x,y) + m(s,y,z) \geq m(t+s,x,z).
\]
\item[(b)] $m$ is periodic, that is, $x,y\in \R^n$, $w\in \Z^n$, and $t>0$,
\[
m(t,x+w,y+w)=m(t,x,y).
\]
\item[(c)] For $t>0$, and $|y| \leq Ct$,
\[
m(2t,0,2y)  \leq 2m (t,0,y)+C.
\]
\item[(d)] For $t>0$, and $|y| \leq Ct$,
\begin{equation}\label{1bdd}
\ol m (t,0,y) \leq m(t,0,y)+C.
\end{equation}
\end{itemize}
\end{lem}
The proof of this lemma is standard and is hence omitted.
See \cite{Tran} for details.
We next show that $m$ is essentially superadditive, which is an extension of \cite[Lemma 4]{Burago}. The subtle point in Cauchy problem is that we need to keep the same amount of time $t$ when joining different pieces of curves. 

\begin{lem}\label{lem:2bdd}
Assume \eqref{quad-growth}--\eqref{quad-growth-L}.
Then, for $t>n$ and $y\in \R^n$ with $|y|\leq Ct$,
\begin{equation}\label{superadd}
2 m(t,0,y) \leq m(2t,0,2y)+C.
\end{equation}
In particular,
\begin{equation}\label{2bdd}
m(t,0,y) \leq \ol m (t,0,y) +C.
\end{equation}
Here, $C>0$ is a universal constant depending only on $L$ and $n$.
\end{lem}

\begin{proof}
It is enough to prove \eqref{superadd}. 
Hereafter, $C>0$ represents a universal constant depending only on $L$ and $n$. 
By considering $\alpha(s)=sy/t$ for $s\in [0,2t]$, we deduce that $m(2t,0,2y)\leq Ct$.
 Thanks to \cite[Appendix D]{Tran}, there exists $\gam:[0,2t] \to \R^n$ with $\gam(0)=0$, $\gam(2t)=2y$ such that
\begin{equation}\label{eq:prep}
m(2t,0,2y)=\int_0^{2t} L(\gam(s),\dot \gam(s))\,ds \leq Ct.
\end{equation}

\smallskip

Let $\xi(s)=(\gam(s),s)$ for $s \in [0,2t]$.
By Lemma \ref{lem:top}, we are able to find a collection of disjoint intervals $\{[a_i,b_i]\}_{1 \leq i \leq k} \subset [0,2t]$ with $k \leq \frac{n+2}{2}$ such that
\[
\sum_{i=1}^k (\xi(b_i)-\xi(a_i)) = \frac{\xi(2t)-\xi(0)}{2}=(y,t).
\]
Rearranging and shifting $\gam$ on $\{[a_i,b_i]\}_{i=1}^k$ in a periodic way in space to get $\tilde \gam:(0,t) \to \R^n$ such that, for $t_0=0$, $t_j = \sum_{i=1}^j (b_i - a_i)$ for $1\leq j \leq k$,
\begin{itemize}

\item $\tilde \gam(0^+) \in [0,1]^n$;

\item $\tilde \gam|_{(t_{j-1}, t_j)}$ is a periodic shift of $ \gam|_{(a_j, b_j)}$ for $1\leq j \leq k$;

\smallskip

\item for $1\leq j \leq k-1$, $\tilde \gam(t_j^+) - \tilde \gam(t_j^-) \in [0,1]^n$, which gives
\[
\left|\tilde \gam(t_j^+) - \tilde \gam(t_j^-) \right| \leq \sqrt{n}; 
\]

\item 
\[
\sum_{i=1}^k (\tilde \gam(t_i^-) - \tilde \gam(t_{i-1}^+))=y.
\]

\end{itemize}
Set $\tilde \gam(0^-)=0$, and $\tilde \gam(t^+)=y$.
We now use $\tilde \gam$ to create $\eta \in \AC([0,t],\R^n)$ with $\eta(0)=0, \eta(t)=y$, and
\begin{equation}\label{eq:connect}
\int_0^{t} L(\eta(s),\dot \eta(s))\,ds\leq \int_0^{t} L(\tilde \gam(s),\dot {\tilde \gam}(s))\,ds+C.
\end{equation}
We create $\eta$ by using $\tilde \gam$ and connectors.
All connectors are straight lines with constant velocity.

\smallskip

If $t<n$, then we simply let $\eta$  be the connector connecting $0$ to $y$, that is, $\eta(s)=sy/t$ for $s\in [0,t]$.
It is clear that \eqref{eq:connect} holds.

\smallskip

Let us now consider the case where $t\geq n$.
By \eqref{eq:prep}, there exists $d\in \{0,1,\ldots, \lfloor t \rfloor -1\}$ such that
\begin{equation*}
\int_d^{d+1} L(\tilde \gam(s),\dot {\tilde \gam}(s))\,ds \leq C.
\end{equation*}
By the bound \eqref{quad-growth-L}, we yield
\begin{equation}\label{eq:con-1}
\int_d^{d+1} |\dot {\tilde \gam}(s)|^2\,ds \leq C.
\end{equation}
Let us now rescale $\tilde \gam$ on $[d,d+1]$ to save an amount of time $1/2$.
Denote by
\[
\alpha(s)=
\begin{cases}
{\tilde \gam}(s) \qquad &\text{ for } 0\leq s \leq d,\\
{\tilde \gam}(d+2(s-d))\qquad &\text{ for } d\leq s \leq d+\frac{1}{2},\\
{\tilde \gam}(s+\frac{1}{2})  \qquad &\text{ for } d+\frac{1}{2} \leq s \leq t-\frac{1}{2}.
\end{cases}
\]
Thanks to \eqref{eq:con-1}, 
\begin{equation}\label{eq:con-2}
\left| \int_0^{t} L(\tilde \gam(s),\dot {\tilde \gam}(s))\,ds - \int_0^{t-\frac{1}{2}} L(\alpha(s),\dot \alpha(s))\,ds \right| \leq C.
\end{equation}
We next create $k+1$ connectors, each takes an amount of time $1/(2k+2)$ connecting $\tilde \gam(t_j^-)$ to $\tilde \gam(t_j^+)$ for $0\leq j \leq k$.
We then glue the pieces of $\alpha$ together with these $k+1$ connectors to get the desired path $\eta$.
See Figure \ref{fig1}.
\begin{center}
\includegraphics[scale=0.55]{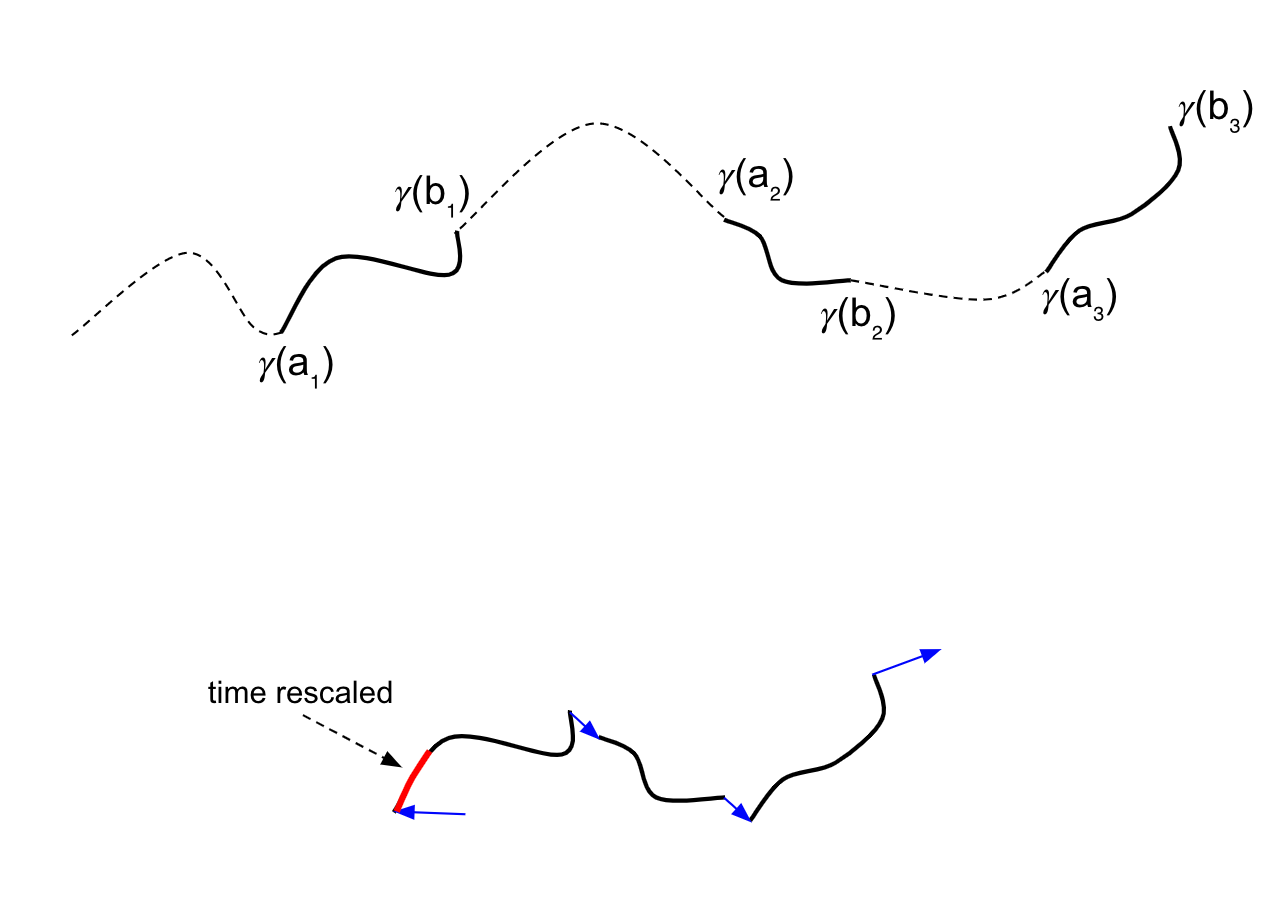}
\captionof{figure}{Formation of the curve $\eta$}\label{fig1}
\end{center}
In light of \eqref{eq:con-2}, \eqref{eq:connect} holds true. 
Combining \eqref{eq:prep} and \eqref{eq:connect}, we arrive at
\[
2 m(t,0,y) \leq m(2t,0,2y)+C.
\]

%We only explain in the following the first step of dealing with $\tilde \gam|_{[t_0,t_1)}$.
%If $\tilde \gam(t_0^+)=0$, then there is nothing to do, and we just let $\eta|_{[t_0,t_1)}=\tilde \gam|_{[t_0,t_1)}$.
%We now assume $\tilde \gam(t_0^+) \neq 0$.
%If $t_1-t_0 <1$, then we replace $\tilde \gam|_{[t_0,t_1)}$ by a connector connecting $0$ to $\tilde \gam(t_1^-)-\tilde \gam(t_0^+)$ to create $\eta|_{[t_0,t_1)}$.
%If $t_1 - t_0 \geq 1$, then we replace $\tilde \gam|_{[0,1/3]}$ by a connector connecting $0$ to $\tilde \gam(1/3)$ to create $\eta|_{[0,1/3]}$.
%Next, let $\eta|_{[1/3,t_1-1/3]}=\tilde \gam|_{[1/3,t_1-1/3]}$.
%Finally, we replace $\tilde \gam|_{[t_1-1/3,t_1)}$ by a connector connecting $\tilde \gam(t_1-1/3)$ to $\tilde \gam(t_1^-)-\tilde \gam(t_0^+)$ to create $\eta|_{[t_1-1/3,t_1)}$.
%Since  the number of connectors created is at most $n+2$

\end{proof}

\begin{proof}[Proof of  Theorem \ref{thm:main}]
We combine \eqref{1bdd} and \eqref{2bdd} to yield, for $s>0$ and $y\in \R^n$ with $|y| \leq Cs$,
\begin{equation}\label{3bdd}
|m(s,0,y)-\ol m(s,0,y)| \leq C.
\end{equation}
By scaling and translation, it suffices to obtain the result for $(x,t)=(0,1)$.
As $g \in \Lip(\R^n)$,
\begin{equation}\label{g-bound}
|g(x)|\leq |g(x)-g(0)| + |g(0)| \leq C(|x|+1) \quad \text{ for all } x\in \R^n.
\end{equation}
The optimal control formula \eqref{oc} gives us that
\[
u^{\ep}(0,1)=\inf_{\substack{  \eta(0)=0 \\ \eta \in \AC([-\ep^{-1},0])}}\left\{g\left(\ep\eta \left(-\ep^{-1}\right)\right)+\ep\int_{-\ep^{-1}}^{0}L(\eta(t), \dot \eta(t))\,dt\right\}.
\]
Thanks to \eqref{quad-growth-L} and Jensen's inequality,
\[
\ep\int_{-\ep^{-1}}^{0}L(\eta(t), \dot \eta(t))\,dt \geq  \ep\int_{-\ep^{-1}}^{0}\left(\frac{|\dot \eta(t)|^2}{2} - K_0 \right)\,dt \geq  {1\over 2}\ep^2\left|\eta (-\ep^{-1})\right|^2-K_0.
\]
Combining the above with \eqref{g-bound},  we obtain that, for  $C>0$ depending only on $L$, $\|Dg\|_{L^\infty(\R^n)}$, and $n$,
\begin{align*}
u^{\ep}(0,1)&=\inf_{\substack{\eta(0)=0, \\ \ep\left|\eta (-\ep^{-1})\right|\leq C}}\left\{g\left(\ep\eta \left(-\ep^{-1}\right)\right)+\ep\int_{-\ep^{-1}}^{0}L(\eta(t),\dot \eta(t))\,dt\right\}\\
&=\inf_{|y|\leq C} (g(y)+\ep m(\ep^{-1},\ep^{-1}y,0))\\
&=\inf_{|y|\leq C} (g(y)+\ol m(1,0,-y)) + O(\ep)\\
&= u(0,1)+ O(\ep).
\end{align*}
We used \eqref{3bdd} in the second last equality.
The proof is complete.
\end{proof}

\begin{rem}
The rate of convergence $O(\ep)$ is indeed optimal.
See \cite[Proposition 4.3]{MTY} for an example demonstrating this.  
This  method from \cite{Burago} is quite robust.  
By necessary modifications, it is not hard to see that the optimal convergence rate $O(\ep)$  can also be established for  homogenizable noncoercive Hamiltonians (e.g., the G-equation from the modeling of turbulent combustion  is one of the most physically relevant examples  (see \cite{CNS,Pet00,XY1})).  
This is left to interested readers as an exercise.   
The general principle is, in the convex setting, as long as the homogenization holds, the convergence rate is likely $O(\ep)$ by properly adapting the methods in \cite{Burago}. 
See also \cite{JTY,Tu}. 
\end{rem}

\section{Proof of Theorem \ref{thm:main2}}\label{sec:proof2}

We assume the setting of Theorem \ref{thm:main2}.
Fix $(x,t)=(0,1)$.
Since ${\rm osc}(V) \geq 1$, according to \cite{ATY1},  the effective Hamiltonian is the same as the effective Hamiltonian associated with 
\[
\tilde H(y,p)=\max\{0, |p|-1\}+V(y).
\]

\subsection{Upper bound}
The upper bound follows easily from the comparison principle and Theorem \ref{thm:main}. 
\begin{lem}
There exists $C>0$ depending only on $\|Dg\|_{L^\infty(\R^n)}$ and $V$ such that
\[
u^\ep(0,1) \leq u(0,1)+ C\ep.
\]
\end{lem}

\begin{proof}
Let $\tilde u^\ep$ be the solution to
\[
\begin{cases}
\tilde u_{t}^{\ep}+\tilde H\left({x\over \ep},\tilde u_{x}^{\ep}\right)=0 \quad &\text{ in } \R^n \times (0,\infty),\\
\tilde u^{\ep}(x,0)=g(x) \quad &\text{ on } \R^n.
\end{cases}
\]
By Theorem \ref{thm:main}, we get that
\[
|\tilde u^\ep(0,1) - u(0,1)| \leq C\ep.
\]
On the other hand, as $H \geq \tilde H$, we imply $\tilde u^\ep \geq u^\ep$.
Thus,
\[
u^\ep(0,1) \leq u(0,1)+ C\ep.
\]
\end{proof}

\subsection{Lower bound} We use game theory representation formulas to derive the lower bound. 
Since
\[
H(x,p)=\max\{|p|-1,1-|p|\}+V(x),
\]
we imply
\[
K(x,p)=-H(x,p)=\min_{\substack{i=1,2 \\ |b| \leq 1}} \, \max_{|a| \le 1} \left( f(a,(i, b))\cdot p + h(i,x)\right).
\]
Here,
\[
f(a, (1,b))=a, \qquad f(a,(2, b))=b,
\]
and
\[
h(1,x)=-1-V(x), \qquad h(2,x)=1-V(x).
\]
The control sets for players I and II are $A=\ol B_1(0)$, and $\cB=\{1,2\}\times \ol B_1(0)$, respectively.
Then, due the representation formula in \cite{ES, Tran}, 
\[
u^{\ep}(x,t)= \sup_{\al \in \Sigma_0} \inf_{\beta \in \cB_t} \left(g(\xi(t))+\int_{0}^{t}h\left(i,{\xi(s)\over \ep}\right)\,ds\right).
\]
Here, $\beta(s)=(i,b(s)) \in \{1,2\}\times \ol B_1(0)$, and
\begin{align*}
\cA_t &= \left\{a:[0,t] \to \ol B_1(0)\,:\, a \text{ is measurable} \right\},\\
\cB_t &=\left\{\beta:[0,t] \to \{1,2\}\times \ol B_1(0)\,:\, \beta \text{ is measurable}\right\},\\
\Sigma_0 &=\left\{ \alpha: \cB_t \to \cA_t \text{ nonanticipating}\right\}.
\end{align*}
In the above, nonanticipating means that, for all $\beta_1(\cdot),\beta_2(\cdot) \in \cB_t$ and $s\in (0,t)$,
\[
\beta_1(\cdot) = \beta_2(\cdot) \text{ on } [0,s) \quad \Longrightarrow \quad \alpha[\beta_1](\cdot) = \alpha[\beta_2](\cdot) \text{ on } [0,s).
\]
\begin{lem}
There exists $C>0$ depending only on $\|Dg\|_{L^\infty(\R^n)}$ and $V$ such that
\[
u^\ep(0,1) \geq u(0,1) - C\ep.
\]
\end{lem}

\begin{proof}
Pick $y\in [0,1]^n$ such that $V(y)=\min_{\T^n} V \leq -1$.
Note first that
\[
u^\ep(0,1) \geq u^\ep(\ep y,1) - C\ep.
\]
We have
\[
u^{\ep}(\ep y,1)= \sup_{\al \in \Sigma_0} \inf_{\beta \in \cB_1} \left(g(\xi(1))+\int_{0}^{1}h\left(i,{\xi(s)\over \ep}\right)\,ds\right),
\]
where
\[
\begin{cases}
\dot \xi(s)=f(\alpha[\beta](s),\beta(s)),\\
\xi(0)=\ep y.
\end{cases}
\]
Player I can play the game only if Player II decides to choose $i=1$.
We design a specific strategy for Player I as following, which can be reduced to a control path associated with $\tilde H$. 
Let us say that the game is switched to Player I at time $s\in [0,1]$.
\begin{itemize}
\item[Step 1.] If the current position of the game 
\[
\xi(s) = \ep y,
\] 
then Player I stays put, that is, $\alpha[\beta](r)=0$ for $r\geq s$ until the game switches back to Player II.

\item[Step 2.] If $\xi(s) \neq \ep y$ then it means that $\xi(s) \neq \xi(0)$, then Player I traces back the step of Player II until the position of the game reaches $\ep y$ or the game is switched back to Player II.
More precisely, 
\[
\xi(s+r)=\xi(s-r) \quad \text{ for } r>0,
\]
and this happens until either Player II choose $i=2$ or $\xi(s+r)=\ep y$ for the smallest $r>0$.
If the latter happens, go to Step 1.

\item[Step 3.] If Step 2 stops at $r_0>0$, then delete $\xi$ on the interval $[s-r_0,s+r_0]$.
If $\xi$ stays put at $\ep y$, then remove this stay put piece as well.
\end{itemize}
We have two key observations.
Firstly, if $i=2$ and $\xi(s) = \ep y$ for $s\in [t_1,t_2]$, then
\[
\int_{t_1}^{t_2}h\left(2,{\xi(s)\over \ep}\right)\,ds = \int_{t_1}^{t_2} (-1-V(y))\,ds \geq 0.
\]
Secondly, if Player II plays $\xi(s-r)$ for $r\in [0,r_0]$, and Player I plays $\xi(s+r)=\xi(s-r)$ for $r\in [0,r_0]$, then
\[
\int_{-r_0}^{0}\left(1-V\left({\xi(s)\over \ep}\right)\right)\,ds + \int_{0}^{r_0} \left(-1-V\left({\xi(s)\over \ep}\right)\right)\,ds  \geq 0.
\]
We can thus safely remove $\xi$ on the interval $[s-r_0,s+r_0]$ without losing anything.
Also, if Step 1 above happens, we remove this stay put path from $\xi$ as well.
By doing it this way, we have removed all the moves of Player I in the game.
Let $\eta:[0,T]\to \R^n$ be the remaining path of $\xi$ after removing all these repeated pieces for some $T\leq 1$.
Then,
\begin{align*}
u^{\ep}(\ep y,1) &\geq \inf_{|\dot\eta(s)|\leq 1} \left(g(\eta(T))+\int_{0}^{T}\left(1-V\left({\eta(s)\over \ep}\right)\right)\,ds\right)\\
&\geq u(\ep y, T) - C\ep \geq u(0,1) - C\ep.
\end{align*}
In the last inequality above, we used the fact that $u_t \leq 0$, and so $u(\cdot,T) \geq u(\cdot,1)$.
\end{proof}

\begin{rem}\label{smallosc} 
For the small oscillation case ${\rm osc}(V) \geq 1$,  if we can show that the number of switches between two players is uniformly bounded, the convergence rate should also be $O(\ep)$. 
This will be investigated in a future project. 
For nonconvex Hamilton-Jacobi equations, it seems  that identifying optimal convergence rates might consist of  (1)   understanding the shape of $\ol H$ and (2)  designing a corresponding game strategy.  
When $n=1$,  the shape of $\ol H$ has been well characterized in \cite{ATY2} for different purposes. 
See \cite{QTY} for some special cases in multi dimensions.   
Due to  different types of nonconvexity, it might be a good idea to look at examples that have physical meaning in the context of homogenization.  One such example is  the G-equation model when strain effect is included (see \cite{Pet00,XY2}), where the Hamiltonian is qualitatively  similar to  $H(x,p)=|p|+A{\cos 2\pi x_1\cos 2\pi x_2(|p_1|-|p_2|)}$ for a given constant $A$ (flow intensity) when $n=2$. 

\end{rem}

%%%%%%%%%%%%%%%%%%%%%%%%%%%%%%%%%%%%%%%%%%%%%%%%%%%%%%%%%%%%%%

%%%%%%%%%%%%%%%%%%%%%%%%%%%%%%%%%%%%%%%%%%%%%%%%%%%%%%%%%%%%%%%%%%%%%%%%

\end{document}